\documentclass[10pt,journal]{IEEEtran}
\usepackage{graphicx}          % För att få grafik att fungera.
\usepackage{amsmath}           % Innehåller vissa matematiska
\usepackage{amssymb}
\usepackage{algorithm}
\usepackage{url}
\usepackage{algpseudocode}
\usepackage[latin1]{inputenc}  % Ser till så att svenska tecken
\usepackage{subfigure}
\usepackage{cite}

\usepackage{amsthm}

\theoremstyle{definition}

\theoremstyle{remark}

\theoremstyle{plain}

\DeclareMathOperator{\E}{E}
\DeclareMathOperator*{\argmin}{arg\,min}

\newcommand{\mbf}[1]{\mathbf{#1}}
\newcommand{\mbs}[1]{\boldsymbol{#1}}

\begin{document}

\title{Self-Localization of Asynchronous Wireless Nodes with Parameter Uncertainties} %Alternative?
\author{Dave Zachariah, Alessio De Angelis, Satyam Dwivedi and Peter Händel\thanks{The authors are with the ACCESS Linnaeus Centre, KTH Royal Institute of
Technology, Stockholm. E-mail: \{davez,  ales, dwivedi, ph\}@kth.se.
Parts of this work have been funded by The Swedish Agency for Innovation Systems (VINNOVA). Copyright (c) 2012 IEEE. Personal use of this material is permitted. However, permission to use this material for any other purposes must be obtained from the IEEE by sending a request to pubs-permissions@ieee.org.}}

%\date{2000-10-10} % Om du t.e.x. vill skriva in datum själv.

\maketitle

\begin{abstract}
We investigate a wireless network localization scenario in which the need for synchronized nodes is avoided.
It consists of a set of fixed anchor nodes transmitting according to a given sequence and a self-localizing receiver node.
The setup can accommodate additional nodes with unknown positions participating in the sequence.
We propose a localization method which is robust with respect to uncertainty of the anchor positions and other system parameters.
Further, we investigate the Cram\'er-Rao bound for the considered problem and show through numerical simulations that the proposed method attains the bound.
\end{abstract}


\section{Introduction}

The rapid emergence of wireless sensor network (WSN) applications, often requiring accurate knowledge of the node positions, has created remarkable research interest for the WSN localization field \cite{PatwariEtAl2005}. In such applications, nodes need to estimate their own position, i.e., perform \emph{self-localization}, by processing noisy range measurements with respect to a set of anchors in a decentralized fashion.

Typically, the anchor nodes positions are assumed to be known exactly, thus neglecting a potential source of error in many application scenarios.
In the literature, the problem of sensor network localization in the presence of anchor position uncertainty has been posed as a maximum likelihood (ML) estimation problem, cf. \cite{LuiEtAl2009}. However, due to the highly nonlinear likelihood function, a closed-form expression of the ML estimator cannot be derived. Therefore the solution has been approached using semi-definite programming in \cite{LuiEtAl2009} and second-order cone programming in \cite{ShiraziEtAl2011}.
Further, in \cite{ZhengAndWu2010}, the authors studied the joint network localization and time synchronization problem with inaccurate anchors, pointing out the strong interdependence between the timing and positioning aspects. All the above mentioned works assume that the timing measurement noise variance is known.

In this letter, we consider the problem of self-localization in sensor networks with uncertainty in anchor position, without assuming knowledge of the timing noise level. More importantly, we propose a system configuration which allows for a receiver node to perform self-localization even without time synchronization of the anchor nodes. In such a configuration the transceiving anchor nodes transmit signals in a predetermined sequence, one after the other, with a turn-around delay which we assume is not perfectly known at the receiver. By passively listening to the transmissions and exploiting available prior knowledge, the receiver can localize itself as well as the transceivers that participate in the sequence. This setup can accommodate the self-localization of an indefinite number of receiver nodes.

We propose an iterative maximum a posteriori (MAP) estimator that effectively copes with uncertainty in the deployment of anchors and with tolerance of the turn-around delay of the asynchronous transceivers. It also provides delay estimates which may be useful for hardware calibration.
Moreover, we study the fundamental performance bounds for the problem considered by deriving the hybrid Cram\'er-Rao bound.
Finally, we evaluate the performance of the proposed estimator by numerical simulations of an ultra-wideband sensor network setup as an application example \cite{DeAngelisEtAl2013}.

\emph{Notation:} $\| \mbf{x} \|_{\mbf{W}} = \sqrt{\mbf{x}^\top \mbf{W} \mbf{x}}$ is the weighted norm, where $\mbf{W}$ is positive definite.% $\mbf{A} \oplus \mbf{B}$ denotes the direct sum of matrices $\mbf{A}$ and $\mbf{B}$.

\section{Problem formulation}

We consider a wireless network of $N-1$ transceiving nodes operating
asynchronously with local clocks. They transmit signals
according to a known sequence, denoted $\mathcal{T}$, set across the
network \cite{DwivediEtAl2012}; when the next node in the sequence receives a signal, it
transmits in return after a certain delay. On this basis, the goal is
to achieve self-localization of the $N$th node, which is a
passive receiver that knows $\mathcal{T}$.

The signals are assumed to have a resolvable temporal signature that
allows for timing events, e.g., pulses, symbol boundaries, etc. and
the propagation velocity $c$ is known. Let $\mbf{x}_i \in
\mathbb{R}^d$ denote the position of node $i$, where $d=2$ or 3, and
$\rho_{i,j} \triangleq \| \mbf{x}_i - \mbf{x}_j \|_2$ denote the range
between nodes $i$ and $j$. Then the observed time interval between a
\emph{pair} of signals received at node $N$, involving transceiving
nodes $i$ and $j$, is modeled by
\begin{equation}
\begin{split}
y^{(i,j)} &= \frac{1}{c}\rho_{i,j} + \delta_j + \frac{1}{c}\rho_{j,N} - \frac{1}{c}\rho_{i,N} + w^{(i,j)},
\end{split}
\label{eq:observationmodel}
\end{equation}
where $\delta_j$ denotes the turn-around delay at node $j$,
cf. \cite{DwivediEtAl2012, GholamiEtAl2012}. The delay $\delta_j$ is generated without any common time reference between the nodes, and is therefore asynchronous. The noise $w^{(i,j)}$
arises from three uncorrelated timing measurements, one at node $j$ and two at node $N$, and is modeled as zero-mean Gaussian with unknown variance $\E[(w^{(i,j)} )^2] = \sigma^2$. The next observed time interval
\begin{equation}
\begin{split}
y^{(j,k)} &= \frac{1}{c}\rho_{j,k} + \delta_k + \frac{1}{c}\rho_{k,N} - \frac{1}{c}\rho_{j,N} + w^{(j,k)},
\end{split}
\end{equation}
uses one timing measurement from the previous observation. Modeling the timing noise variance across nodes equally we have the correlation $\E[w^{(i,j)}  w^{(j,k)}] = \sigma^2 / 3$ for all consecutive observations.

Prior knowledge of the positions will be modeled as $\mbf{x}_i \sim
\mathcal{N}(\mbs{\mu}_i, \mbf{P}_i )$, where $\mbs{\mu}_i$ and $\mbf{P}_i$ are known, $\forall i$. Setting $\mbf{P}^{-1}_i = \mbf{0}$,
leads to a noninformative prior, $p(\mbf{x}_i) \propto 1$ \cite{Tiao&Zellner1964}. Further, to avoid signal collisions
it is necessary that the delays exceed $\rho_{\text{max}}/ c$, where
$\rho_{\text{max}}$ is the maximum range between any pair of
transceivers and can easily be ensured in any bounded localization
scenario. Whilst the delay in each transceiver node may be set to some
nominal value $\mu_\delta$, the actual delay $\delta_j$ will deviate
due to hardware imperfections. We model this as $\delta_j
\sim \mathcal{N}(\mu_\delta, \sigma^2_{\delta})$, with a given
$\sigma^2_{\delta}$. For the unknown noise variance, we assume a noninformative prior $p(\sigma^2) \propto 1/\sigma^2$ \cite{Tiao&Zellner1964}.

For notational simplicity we write $\mbs{\theta} \triangleq
[\mbf{x}^\top_1 \cdots \mbf{x}^\top_N]^\top \in \mathbb{R}^{dN}$ and
$\mbs{\delta} \triangleq [\delta_1 \cdots \delta_{N-1}]^\top \in
\mathbb{R}^{N-1}$. The goal is to estimate $\mbs{\theta}$,
$\mbs{\delta}$ and $\sigma^2$ from a set of $M$ observations $\{ y_m \}^M_{m=1}$.

\section{MAP estimator} \label{sec:MAP}

For a given sequence $\mathcal{T}$, the pair of nodes involved in each
observed time interval $\{ y_m \}^M_{m=1}$ is known. In vector form,
the observation model is
\begin{equation}
\mbf{y} = c^{-1} \mbf{H} \mbf{g}(\mbs{\vartheta}) + \mbf{w} \in \mathbb{R}^M,
\end{equation}
where $\mbs{\vartheta} \triangleq [\mbs{\theta}^\top \;
\mbs{\delta}^\top]^\top \in \mathbb{R}^T$, and $T = dN + N-1$. The noise follows $\mbf{w} \sim \mathcal{N}(\mbf{0}, \sigma^2 \mbf{Q})$, where $[\mbf{Q}]_{i,i} = 1$, $[\mbf{Q}]_{i,j} = \frac{1}{3}$, $\forall i,j$ such that $|i-j| = 1$ and $[\mbf{Q}]_{i,j} = 0$ otherwise.
%\begin{equation*}
%\mbf{Q} =
%\begin{bmatrix}
%1           & \frac{1}{3} &  \\
%\frac{1}{3} & 1           & \ddots \\
% & \ddots   &   \\
% &          &  \\
%\end{bmatrix}
%\end{equation*}
The
nonlinear mapping $\mbf{g}(\mbs{\vartheta}) =
[\mbs{\rho}^\top(\mbs{\theta}) \quad \mbs{\delta}^\top]^\top \in
\mathbb{R}^{N(N-1)/2 + N-1}$ contains all
the pairwise ranges and delays, and $\mbf{H}$ is determined by the
transmission sequence $\mathcal{T}$, cf. \eqref{eq:observationmodel}.
We aim to find the maximum a posteriori (MAP) estimator, i.e., the maximizer of $p(\mbs{\vartheta}, \sigma^2 | \mbf{y})$.

\subsection{Concentrated cost function}
Using Bayes' Rule, the MAP estimate can be computed by maximizing $J(\mbs{\vartheta}, \sigma^2) = \ln p( \mbf{y}| \mbs{\vartheta}, \sigma^2 ) +  \ln p(\sigma^2)  + \ln p( \mbs{\vartheta} )$. Further, define
\begin{equation}
\begin{split}
J_1(\mbs{\vartheta}, \sigma^2) &\triangleq \ln p( \mbf{y}| \mbs{\vartheta}, \sigma^2 ) +  \ln p(\sigma^2) \\
&= -\frac{M+2}{2} \ln \sigma^2 - \frac{1}{2 \sigma^2} \| \mbf{y} - c^{-1}
\mbf{H} \mbf{g}(\mbs{\vartheta}) \|^2_{\mbf{Q}^{-1}} + K_1,
\end{split}
\label{eq:J_1}
\end{equation}
where $K_1$ is a constant. Similarly,
\begin{equation}
\begin{split}
J_2(\mbs{\vartheta}) &\triangleq \ln p(  \mbs{\vartheta} ) = -\frac{1}{2} \| \mbs{\vartheta} - \mbs{\mu} \|^2_{\mbf{P}^{-1}} + K_2,
\end{split}
\label{eq:J_2}
\end{equation}
where $K_2$ is a constant. Maximizing \eqref{eq:J_1} with respect to $\sigma^2$ yields the
estimate $\hat{\sigma}^2 = \| \mbf{y} - c^{-1} \mbf{H}
\mbf{g}(\mbs{\vartheta}) \|^2_{\mbf{Q}^{-1}} / ( M+2)$. Inserting this back into
\eqref{eq:J_1}, and combining with \eqref{eq:J_2},
results in a concentrated cost function. The MAP estimator is then given by
\begin{equation}
\hat{\mbs{\vartheta}} = \argmin_{\mbs{\vartheta} \in \mathbb{R}^T}
V(\mbs{\vartheta})
\label{eq:V_map},
\end{equation}
where
\begin{equation} \label{eq:cost_func}
V(\mbs{\vartheta}) \triangleq \frac{1}{2}  \ln \| \mbf{y} - c^{-1} \mbf{H} \mbf{g}(\mbs{\vartheta}) \|^2_{\mbf{Q}^{-1}} + \frac{\beta}{2} \| \mbs{\mu} - \mbs{\vartheta} \|^2_{\mbf{P}^{-1}},
\end{equation}
$\beta=1/(M+2)$, $\mbs{\mu} = [\mbs{\mu}^\top_1 \: \cdots \: \mbs{\mu}^\top_N \;
\mu_{\delta} \mbf{1}^{\top}_{N-1}]^\top$ and $\mbf{P}^{-1} = \text{diag}(
\mbf{P}^{-1}_1 , \cdots , \mbf{P}^{-1}_N ,
\sigma^{-2}_{\delta} \mbf{I}_{N-1})$.

\subsection{Iterative solution}

The optimization problem \eqref{eq:V_map} does not lend
itself to a closed-form solution. A standard method for solving \eqref{eq:V_map} is gradient descent. However, its performance is heavily dependent on the user-defined step length. Instead we propose an iterative
solution by first exploiting the linearization around an initial
estimate $\hat{\mbs{\vartheta}}_\ell$, i.e., $\mbf{g}(\mbs{\vartheta})
\simeq \mbf{g}(\hat{\mbs{\vartheta}}_\ell) +
\mbs{\Gamma}(\hat{\mbs{\vartheta}}_\ell) \tilde{\mbs{\vartheta}} $,
where $\tilde{\mbs{\vartheta}} \triangleq \mbs{\vartheta} -
\hat{\mbs{\vartheta}}_{\ell}$ is the iteration increment and $\mbs{\Gamma}(\mbs{\vartheta})$ is
the Jacobian of $\mbf{g}({\mbs{\vartheta}})$. Then we may write a
cost function which approximates \eqref{eq:cost_func} as
\begin{equation}
\begin{split}
V_\ell(\tilde{\mbs{\vartheta}} ) & \triangleq %\ln \| \mbf{y} - c^{-1} \mbf{H}\mbf{g}_\ell - c^{-1} \mbf{H} \mbs{\Gamma}_\ell \tilde{\mbs{\vartheta}}  \|^2_2 + \beta \| \mbs{\mu} - \hat{\mbs{\vartheta}}_\ell - \tilde{\mbs{\vartheta}} \|^2_{\mbf{P}^{-1}} \\
\frac{1}{2} \ln \| \tilde{\mbf{y}}_\ell - \mbf{G}_\ell \tilde{\mbs{\vartheta}}   \|^2_{\mbf{Q}^{-1}} + \frac{\beta}{2} \| \tilde{\mbs{\mu}}_\ell - \tilde{\mbs{\vartheta}} \|^2_{\mbf{P}^{-1}},
\end{split}
\label{eq:V_map_approx}
\end{equation}
where we introduce $\mbf{g}_\ell = \mbf{g}(\hat{\mbs{\vartheta}}_\ell)$,
$\mbs{\Gamma}_\ell = \mbs{\Gamma}(\hat{\mbs{\vartheta}}_\ell)$,  $\tilde{\mbf{y}}_\ell = \mbf{y} - c^{-1}
\mbf{H}\mbf{g}_\ell$, $\mbf{G}_\ell  = c^{-1} \mbf{H}
\mbs{\Gamma}_\ell$ and $\tilde{\mbs{\mu}}_\ell  = \mbs{\mu} - \hat{\mbs{\vartheta}}_\ell$.
We aim at iteratively updating the estimate $\hat{\mbs{\vartheta}}_{\ell}$ by finding the optimal increment $\tilde{\mbs{\vartheta}}$.

To this end, we use a fixed-point iteration. First, we note that the gradient of
$V_\ell(\tilde{\mbs{\vartheta}})$ equals
\begin{equation*}
\partial_\vartheta V_\ell =  \alpha(\tilde{\mbs{\vartheta}}) \mbf{G}^\top_\ell\mbf{Q}^{-1}
( \mbf{G}_\ell \tilde{\mbs{\vartheta}}  -
  \tilde{\mbf{y}}_\ell )  + \beta \mbf{P}^{-1} (
  \tilde{\mbs{\vartheta}} - \tilde{\mbs{\mu}}_\ell ),
\end{equation*}
where $\alpha(\tilde{\mbs{\vartheta}}) \triangleq 1 / \|
\tilde{\mbf{y}}_\ell - \mbf{G}_\ell \tilde{\mbs{\vartheta}}
\|^2_{\mbf{Q}^{-1}}$.
Then, we hold $\alpha(\tilde{\mbs{\vartheta}})$ fixed and solve $ \partial_\vartheta V_\ell = \mbf{0}$, thus obtaining the following fixed-point iteration
\begin{equation}
\begin{split}
\tilde{\mbs{\vartheta}} &:=
\left(\alpha(\tilde{\mbs{\vartheta}})\mbf{G}^\top_\ell \mbf{Q}^{-1}  \mbf{G}_\ell +
  \beta \mbf{P}^{-1} \right)^{-1} \\
&\quad \times (  \alpha(\tilde{\mbs{\vartheta}})
\mbf{G}^\top_\ell \mbf{Q}^{-1} \tilde{\mbf{y}}_\ell + \beta \mbf{P}^{-1}
\tilde{\mbs{\mu}}_\ell  ).
\end{split}
\label{eq:iteratevartheta}
\end{equation}
By iteratively applying \eqref{eq:iteratevartheta}, starting with the zero increment $\tilde{\mbs{\vartheta}}= \mbf{0}$, we converge to a stationary point. The analytical convergence properties are difficult to derive and are beyond the scope of this letter. However, in Section~\ref{sec:results}, we show the convergence properties by numerical evaluation in a practical scenario. Further, we provide a practical method for obtaining an initial estimate $\hat{\mbs{\vartheta}}_0$.
The iterative estimator is summarized in Algorithm~\ref{alg:MAP}, where $\varepsilon$ is the convergence threshold.
\begin{algorithm}
\caption{Iterative MAP estimator} \label{alg:MAP}
\begin{algorithmic}[1]
\State Input: $\mbf{y}$, $\mbs{\mu}$, $\mbf{P}$, $c$ and
$\hat{\mbs{\vartheta}}_0$
\State Set $\ell := 0$ and $\beta = 1/(M+2)$
\Repeat
%\State Form $\mbf{g}(\hat{\mbs{\vartheta}}_\ell)$ and
%$\mbs{\Gamma}(\hat{\mbs{\vartheta}}_\ell)$
\State  $\tilde{\mbf{y}}_\ell = \mbf{y} - c^{-1} \mbf{H}\mbf{g}(\hat{\mbs{\vartheta}}_\ell)$
\State $\mbf{G}_\ell = c^{-1} \mbf{H} \mbs{\Gamma}(\hat{\mbs{\vartheta}}_\ell)$
\State $\tilde{\mbs{\mu}}_\ell  = \mbs{\mu} -
\hat{\mbs{\vartheta}}_\ell$
\State $\tilde{\mbs{\vartheta}} := \mbf{0}$
\State Repeat \eqref{eq:iteratevartheta} until convergence
\State $\hat{\mbs{\vartheta}}_{\ell+1} = \tilde{\mbs{\vartheta}} +
\hat{\mbs{\vartheta}}_{\ell}$, $\ell := \ell +1$
\Until{$\| \hat{\mbs{\vartheta}}_{\ell} - \hat{\mbs{\vartheta}}_{\ell-1}  \|_2 < \varepsilon$}
\State Output: $\hat{\mbs{\vartheta}}_\ell$
\end{algorithmic}
\end{algorithm}

\section{Cramér-Rao bound} \label{sec:HCRB}

Let $\mbs{\eta} \triangleq [\mbs{\theta}^\top \; \mbs{\delta}^\top \;
\sigma^2 ]^\top \in \mathbb{R}^{T+1}$ and $\hat{\mbs{\eta}}$ be any estimator that is conditionally unbiased with respect to the deterministic parameters.
Then its mean square error (MSE) matrix is constrained by the hybrid Cramér-Rao bound (HCRB) \cite{vanTrees2002}, $\mbf{C}_{\tilde{\eta}} \succeq \mbf{J}^{-1}_{\eta} $, where $
\mbf{J}_{\eta} = \mbf{J}^D_\eta + \mbf{J}^P_\eta \in \mathbb{R}^{(T+1) \times (T+1)} $.

Here $\mbf{J}^D_\eta = \E_{\bar{\eta}}[ \mbf{J}^D(\mbs{\eta}) ]$ is
the expected Fisher information matrix, where $\bar{\mbs{\eta}}$ denotes the subset of parameters that are modeled as random quantities and
\begin{equation*}
[\mbf{J}^D(\mbs{\eta})  ]_{i,j} %= c^{-2} \sigma^{-2} \frac{ \partial \mbf{g}^{\top}
%}{ \partial \eta_i } \mbf{H}^\top \mbf{H}  \frac{\partial \mbf{g}
%}{ \partial \eta_j } + \sigma^{-2} \text{tr}\left\{ \frac{\partial \sigma^2}{\partial %\eta_i} \frac{\partial \sigma^2}{\partial \eta_j} \mbf{I}_M \right\}
= \frac{ c^{-2}}{\sigma^{2} } \frac{ \partial \mbf{g}^{\top}
}{ \partial \eta_i } \mbf{H}^\top \mbf{Q}^{-1} \mbf{H}  \frac{\partial \mbf{g}
}{ \partial \eta_j } + \frac{M}{2 \sigma^{4}}  \frac{\partial \sigma^2}{\partial \eta_i} \frac{\partial \sigma^2}{\partial \eta_j} ,
\end{equation*}
as given in \cite{Kay1993}. As the expectation does not have a closed form solution, we evaluate it by Monte Carlo simulation. If a subset of node positions, $\mbs{\theta}_u$, and
the noise level, $\sigma^2$, are treated as deterministic and unknown
parameters, while the remaining parameters, $\mbs{\theta}_a$ and
$\mbs{\delta}$, are random Gaussian then the prior information matrix is given by
\begin{equation*}
\mbf{J}^P_\eta =
\begin{bmatrix}
\mbf{P}^{-1}_{\theta_a} & \mbf{0} & \mbf{0} & \mbf{0} \\
\mbf{0} & \mbf{0} & \mbf{0} & \mbf{0} \\
\mbf{0} & \mbf{0} & \sigma^{-2}_{\delta} \mbf{I}_{N-1} & \mbf{0} \\
\mbf{0} & \mbf{0} & \mbf{0} & 0
\end{bmatrix} \in \mathbb{R}^{(T+1) \times (T+1)} \,,
\end{equation*}
where $\mbf{P}_{\theta_a}$ is the covariance matrix of $\mbs{\theta}_a$.
This division between deterministic and random parameters may be useful to study practical configurations where we lack prior knowledge on the position of a subset of nodes.

\section{Numerical results} \label{sec:results}
To evaluate the performance of the MAP estimator proposed in Section~\ref{sec:MAP} we provide numerical simulation results in an ultra wideband wireless sensor network scenario \cite{DeAngelisEtAl2013}. In particular, we compare the HCRB derived in Section~\ref{sec:HCRB} with the root mean square error (RMSE) of the MAP estimates.

\subsection{Setup}
We consider a network of $N=6$ nodes consisting of one self-localizing passive receiver, without prior knowledge of its position, and $N_a=4$ anchors with $\mbf{P}_{\theta_a}=\sigma_a^2 \mbf{I}_{2 N_a}$. In addition, we consider one auxiliary node, which is a transceiver participating in the sequence with a noninformative position prior. Therefore, the number of unknown-position nodes is $N_u=2$.
The simulation has been performed using the topology shown in Fig.~\ref{fig:ellipses} and setting $\mu_{\delta}=10^{-6}$~s. In each realization, the positions of the anchors and the delays have been randomly generated according to their prior distributions.
Furthermore, we construct the sequence $\mathcal{T}$ to ensure that all pairwise combinations appear at least once. E.g., $\mathcal{T} = \{1,2,1,3,1,4,1,5,2,3 \dots \}$.
 %therefore the length of $\mathcal{T}$ is $2 \left( N - 1 \right) \left( N - 2 \right)$.

The average RMSE of the position and delay estimates are given by
\begin{align*}
\text{RMSE}_{\xi} & \triangleq \frac{1}{N_{\xi}} \sqrt{ \text{tr}{\left\{ \mbf{C}_{\tilde{\xi}} \right\}}},
\end{align*}
where $\xi$ can be either $\mbs{\theta}_u$, in which case $N_{\xi}=N_u$, or $\mbs{\delta}$, and $N_{\xi}=N-1$. Here $\mbf{C}_{\tilde{\xi}}$ is the MSE matrix of $\xi$. We estimated the RMSE from $10^3$ Monte Carlo iterations.

We initialize the MAP estimator of Algorithm~\ref{alg:MAP} with
$\hat{\mbs{\vartheta}}_0 = \mbs{\mu} = [\mbs{\mu}^\top_1 \: \cdots \: \mbs{\mu}^\top_{N_a} \; \mbf{u} \; \bar{\mbs{\mu}} \;
\mu_{\delta} \mbf{1}^{\top}_{N-1}]^\top \,,$
where $\mbf{u}$ is an arbitrary position initialization of the auxiliary node, with $\mbf{u} \neq \mbs{\mu}_i$. The self-localizing node's initial position estimate $\bar{\mbs{\mu}}$ is the centroid of the anchor positions, i.e. $\bar{\mbs{\mu}} = \frac{1}{N_a} \sum_{i=1}^{N_a} \mbs{\mu}_i$. Further, we set $\varepsilon = 10^{-4}$. The same tolerance is used for Step~8 in Algorithm~\ref{alg:MAP}.

\subsection{Results}

\begin{figure}[t]
\centering
\includegraphics[width = 0.85\columnwidth]{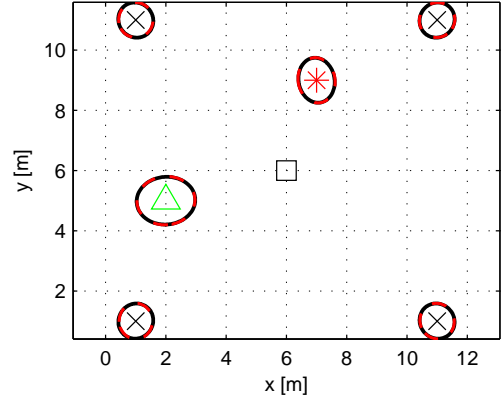}
\caption{Scatter plots of the true node positions and error ellipses for all nodes estimated at the self-localizing node indicated by ($\ast$).
The anchors are denoted by ($\times$), the auxiliary node by ($ \triangle $) and the anchors centroid by ($\square$). The solid black ellipses indicate the CRB and the dashed red ellipses indicate
the MSE performance of the MAP estimator. For visual clarity the sizes
have been scaled to correspond to 99\% confidence ellipses of a zero-mean Gaussian distribution.
Here $ \sigma = 2$~ns, $ \sigma_a = 0.2$~m and $ \sigma_{\delta} = 10$~ns.}
\label{fig:ellipses}
\end{figure}

\begin{figure}[t]
\centering
\includegraphics[width = 0.85\columnwidth]{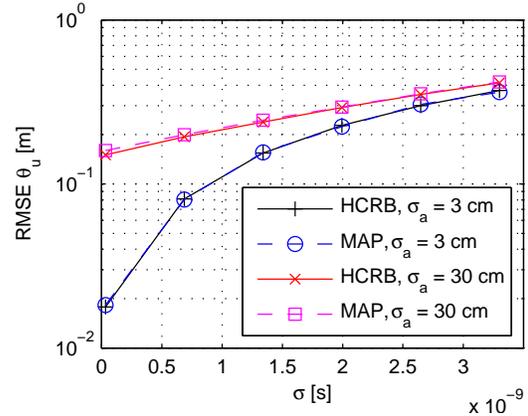}
\caption{Average RMSE of the position estimates of all unknown-position
nodes at the self-localizing node as a function of the standard deviation of the measurement noise $\sigma$. Two values of the prior on the anchors positions are shown. Here, $ \sigma_{\delta} = 10$ ns. }
\label{fig:theta_vs_sigma}
\end{figure}

\begin{figure}[t]
\centering
\includegraphics[width = 0.85\columnwidth]{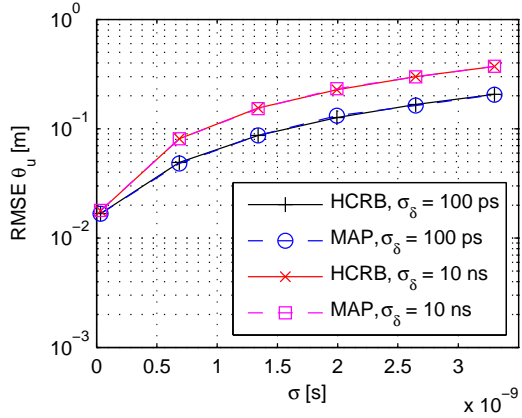}
\caption{Average RMSE of the position estimates of all unknown-position
nodes at the self-localizing node as a function of the standard deviation of the measurement noise $\sigma$. Two values of the prior of the delay are shown. Here, $ \sigma_{a} = 3$ cm. }
\label{fig:theta_vs_sigma_param_delta}
\end{figure}

\begin{figure}[t]
\centering
\includegraphics[width = 0.85\columnwidth]{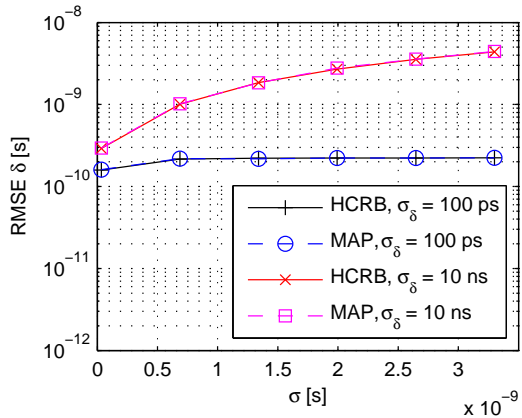}
\caption{Average RMSE of the delay as a function of the standard deviation of the measurement noise $\sigma$. Two values of the prior of the delay are shown. Here, $ \sigma_{a} = 3$ cm.}
\label{fig:delta_vs_sigma_param_delta}
\end{figure}

\begin{figure}[t]
\centering
\includegraphics[width = 0.85\columnwidth]{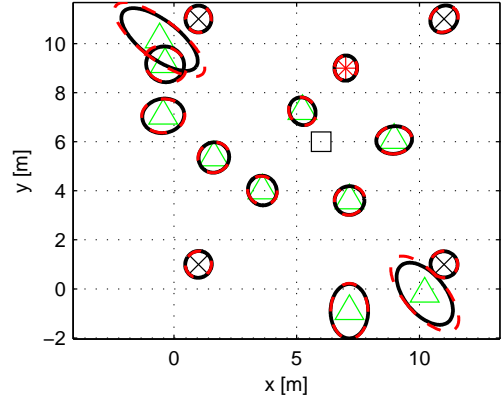}
\caption{Scatter plots of the true node positions and error ellipses for a uniform random placement of 10 auxiliary nodes denoted by ($ \triangle $).
%The solid black ellipses indicate the CRB and the dashed red ellipses indicate
%the MSE performance of the MAP estimator. For visual clarity the sizes
%have been scaled to correspond to 99\% confidence ellipses of a zero-mean Gaussian distribution.
Here $ \sigma = 2$ ns, $ \sigma_a = 0.2$ m and $ \sigma_{\delta} = 10$ ns.}
\label{fig:ellipses_10txrx}
\end{figure}

\begin{figure}[t]
\centering
\includegraphics[width = 0.85\columnwidth]{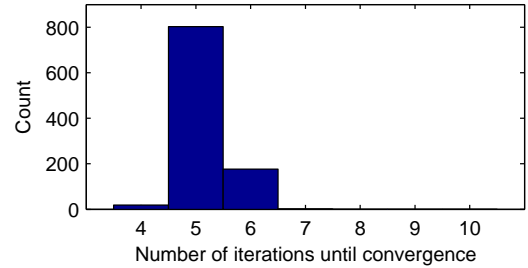}
\caption{Histogram of the number of iterations of the MAP estimator in the scenario in Fig.~\ref{fig:ellipses} for $10^3$ realizations.}
\label{fig:conv}
\end{figure}

Fig.~\ref{fig:ellipses} shows the topology of the network, together with the error ellipses produced by the MAP estimator and given by the HCRB. It can be seen that the proposed estimator attains the bound.
Next, we provide an evaluation of the position estimates as a function of timing noise level $\sigma$. In Fig.~\ref{fig:theta_vs_sigma}, we vary the uncertainty of the anchor positions which is parameterized by $\sigma_a$. In a low noise scenario, we achieve an accuracy of the same order of magnitude as the anchor position prior. In Fig.~\ref{fig:theta_vs_sigma_param_delta}, we vary the uncertainty of the delays set by $\sigma_{\delta}$. The results show that the proposed position-estimation method is robust with respect to deviations from the nominal delay. Similarly, Fig.~\ref{fig:delta_vs_sigma_param_delta} shows the robustness of the delay estimator. We tested a range of deviations that spans two orders of magnitude.
%
%\begin{figure}[t]
%\centering
%\includegraphics[width = 0.85\columnwidth]{fig_1k_delta_vs_sigma_3_30.eps}
%\caption{Average RMSE of the delay as a function of the standard deviation of the timing
%noise. Two values of the prior on the anchors positions are shown. $ \sigma_{\delta} = 10$ ns.}
%\label{fig:delta_vs_sigma}
%\end{figure}
%
Furthermore, the extended simulation setup in Fig.~\ref{fig:ellipses_10txrx} shows that the proposed method can localize $10$ auxiliary nodes which participate in the transmission sequence and are not collocated.
However, the MAP does not attain the HCRB in all scenarios, in particular when nodes are close to anchors as can be seen in Fig.~\ref{fig:ellipses_10txrx}.
Moreover, Fig.~\ref{fig:conv} shows the convergence behavior of the estimator, the average is $5.12$ iterations.

Finally, we repeated the simulations in the scenario in Fig.~\ref{fig:ellipses}, assigning a more conservative prior with respect to the actual variability of the anchor positions. In particular, we defined $\sigma'_a = 10 \sigma_a$, where $\sigma'_a $ is used by Algorithm~\ref{alg:MAP} and $\sigma_a=0.2$~m is used to generate the random anchor positions. We obtained RMSE$_{\mbf{\theta}_u} = 0.37$ m and HCRB$_{\mbf{\theta}_u} = 0.26$ m, indicating inherent robustness to such model mismatches.

%\emph{Reproducible research:} TODO

\section{Conclusion}
In this letter, we have considered a wireless network localization scenario with asynchronous nodes where the entire network topology can be estimated by a self-localizing receiver node.
We have proposed a MAP estimator and compared its performance with the HCRB that we derived for the problem.
The simulation results show that the estimator attains the bound and is robust with respect to uncertainty of the anchor positions, delays and noise level.
In addition to the anchor nodes, we found that the setup allows for the localization of auxiliary nodes participating in the transmission sequence.
\bibliography{refs_self_loc}
\bibliographystyle{ieeetr}

\end{document}